\documentstyle[12pt,amscd,amssymb]{amsart}
\input xypic \xyoption{curve}

\def\centerbmp#1#2#3{\vskip#2\relax\centerline{\hbox to#1{\special
  {bmp:#3 x=#1, y=#2}\hfil}}}

\newcommand{\B}[1]{{\bold#1}} 
\newcommand{\C}[1]{{\mathcal#1}} 

\theoremstyle{plain}

\newtheorem{lem}{Lemma}[section]
\newtheorem{prop}{Proposition}[section]
\newtheorem{conj}{Conjecture}[section]

\theoremstyle{definition}
\newtheorem{defin}{Definition}[section]

\theoremstyle{definition}

\theoremstyle{remark}
\newtheorem{rem}{Remark}[section]

\begin{document}

\title{A Natural Partial Order on The Prime Numbers}         
\author{Lucian M. Ionescu}
\address{Department of Mathematics, Illinois State University, IL 61790-4520}
\email{lmiones@@ilstu.edu}
\date{\today}         

\begin{abstract}
A natural partial order on the set of prime numbers was derived by the author
from the internal symmetries of the primary finite fields \cite{LMI:talk-DMS},
independently of \cite{FKL}, who investigated Pratt trees \cite{Pratt}
used for primality tests.

It leads to a correspondence with the Hopf algebra of rooted trees, 
and as an application, to an alternative approach to the Prime Number Theorem.
\end{abstract}

\maketitle
\tableofcontents

\section{Introduction}       
There is a need for a new insight into the ring of integers 
\cite{YuManin}, p.1, \cite{ShaiHaran}, \cite{PaulaTretkoff}, p.143, etc.
As usual, physics comes to the rescue\footnote{E.g. invariants of knots, 4-manifolds etc.}: 
in a quantum digital universe, the primary modes of vibration (periodic structures)
should be modeled by primary fields $F_p$, reflecting the value of the fine structure constant
via a modern, quantum computing, model of the hydrogen atom \cite{LMI-Rem}.

One way or another, {\em categorification} is the correct approach: instead of investigating 
``sizes'' of objects, study the objects themselves, not their ``shadows'':
$$\diagram
Categories: & \quad Finite\ Sets \ \C{S} \rrto^{\quad Simple\ objects} \dto_{Grothendieck}& &
Spec_\C{S} \dto_{``Tangent\ mapping''}^{(generators)}\\ 
Number \ Systems: & Integers \ \B{Z}\quad \rrto^{Prime\ Numbers} & & \quad Spec(\B{Z}).
\enddiagram$$
Without loss of generality for modeling reality, one can restrict to the study of
the category of finite sets $\C{S}$.
Once a coordinate system is chosen on sets, reducing their symmetry groups,
the category of finite abelian groups is obtained $\C{A}b$.
The main object of study, then, is their symmetries, the {\em internal}
functor $Aut_{\C{A}b}$.

In this article we will introduce a lattice structure on the set of prime numbers,
derived from $Aut_{\C{A}b}$,
leading to a mapping into the Hopf algebra of rooted trees.
It is natural to investigate what are the functions analog to those 
involved in the proof of the Prime Number Theorem, and 
whether on can obtain a better understanding of the latter ``lifting it''
through the canonical POSet morphism from the lattice of primes,
to natural numbers.

In conclusion, the connection between these two primeval structures, primes and 
trees, suggests deeper results in the light of Selberg's trace formula interpretation 
of the Riemann-Mangoldt/Weil exact formula.

Additional explanations regarding the author's approach justifying some of the
above claims will be given in a follow up article
(see also the extended version available on the web).

\section{A partial order on the set of prime numbers}

Primes are sizes of the basic finite fields $p=|F_p|$. These, as abelian groups,
have ``symetries'' which in turn are abelian groups with symmetries etc.
This leads to a sort of a resolution of the basic finite fields.
Its ``shadow'', via size, determines a partial order on the set of primes. 
This partial order, conjecturally a lattice structure, allows to define several gradings, leading to a concept of complexity of a prime number.

\subsection{The Lie-Klein picture of a finite fields}
The automorphisms of $(F_p,+)$ as an abelian group ($\B{Z}$-module) form
a group isomorphic to the corresponding multiplicative group via
the regular representation 
\footnote{Thus these irreducible ``discrete vector spaces''
{\em happen} to be (finite) fields.}.

Its summands are elementary factors of the form $Z_{q^k}$, 
corresponding to the factorization of Euler's totient function value:
$$L:(F_p,\cdot)\to Aut(F_p,+)\cong \prod_{p \in \B{P}} 
\B{Z}_{p^{k(p)}}, \quad \phi(p)=p-1=\prod_{p\in \B{P}} p^{k(p)}.$$
Here a positive integer is determined by its {\em exponent function} $k(p)$,
defined on the set of primes $\B{P}$. 
\footnote{... towards a function field interpretation of the rationals,
and application of duality of algebraic quantum groups.}
For example, $Aut(F_{19})\cong \B{Z}_2\times \B{Z}_{3^2}$, since $19-1=2\cdot 3^2$.

The factorization of $\phi(p)=p-1=2^m \cdot q$ will be called the 
{\em Proth factorization}.

\subsection{Prime extensions and their complexity}
In this way the smaller primes $q$ dividing $p-1$ represent the ``Lie generators''
of the ``space'' $(F_p,+)$, allowing for the definition of a natural partial order.
\begin{defin}
A prime $q$ is called a generator of the prime $p$, denoted $q<<p$ 
iff $q$ divids $p-1$.
\end{defin}
For example $2$ and $3$ are generators of $p=19$.
Fermat primes $p=2^n+1$ have only one generator, $q=2$.

Note that the ``multiplicity'' of a prime in the factorization 
of $p-1$ is disregarded at this stage.
\begin{prop}
The set of prime numbers $(\B{P},<<)$ is a connected partial ordered set (POSet)
with minimal element $2$.
The ``identity'' function:
$$I:(\B{P},<<)\to (\B{N},<), \ I(p)=p$$
is a morphism  of POSets, 
representing a \underline{refinement} of the total ordering of primes by ``size''.
\end{prop}
The properties of this POSet are not clear to the author at this time, e.g.
whether it is a lattice, its cohomology etc..
The ``Euclid's trick'' $N=p \cdot p'+1$ suggests the existence of a $p\cup p'$,
but there seems to be no reason for the uniqueness of $p\cap p'$.

\subsection{The relation with Pratt trees and prime chains}
V. Pratt introduced the concept of prime certificate in \cite{Pratt},
as a collection of data associated to a given number $n$ as his root,
as a proof of primality via the Lucas primility test.
The nodes include primitive roots of $Z/qZ$ besides the divisors $q|p-1$.

The present author arrived to the same notion starting top-down, 
via categorification of primes, investigating the
structure of basic finite fields.

The emphasis in this article is on the global structure of $P$ as a POSet,
and its properties towards a better understanding of multiplicative number theory and
the duality between prime numbers and the non-trivial zeroes of the Riemann zeta function.

The partial order was directly introduced in \cite{FKL}, inspired by Pratt trees.
Chains of primes are paths in Pratt trees.

In brief, disregarding multiplicity in the subtree of partial order determined by $p$ in $P$ coincides with the tree obtained by forgetting the primitive roots in a Pratt tree.

For example the tree $t_{47}$ representing the cutoff of the POSet $P$ at $p=47$
$$\diagram
t(47)= & 47 \dto \rto & 23 \dto \rto & 11 \dto \rto & 5 \dto \rto & 2\\
& 2 & 2 & 2 & 2
\enddiagram$$
is implicit in the Pratt tree (prime certificate or proof) from \cite{Pratt}, p.215.
In the above tree the nodes are labeled with the corresponding primes, although the
``prime value'' of the node can be uniquely determined going up the branches.

The relation between the POSet $P$ and its trees on one hand, 
and Pratt trees/ prime chains on the other, can be used to transport
statistical results regarding their ``abundance'', maximal depth etc.;
for details regarding such results see \cite{FKL}.

The main advantage of thinking globally and introducing the POSet $P$ ``naturally''
(conceptually), 
from its categorical interpretation is the possibility of studying the correlation
of primes, via their common symmetries, from an algebraic point of view 
related to p-adic numbers \cite{LMI-AQG, LMI-rationals, LMI-realsandp-adic}.
This direction of research will be addressed elsewhere.

\subsection{Gradings on the set of primes}
There are a few natural functions defined on $P$, 
including a grading by ``support''.
\begin{defin}
The function $w:P\to \B{N}$, $w(p)=dim(Aut(F_p,+))$ is called 
the {\em grading of primes by support}.
\end{defin}
The term ``support'' is used because of the interpretation of 
the multiplicative (positive) rationals
as functions on $P$: $r=m/n:P\to \B{Z}$,
corresponding to their unique factorization.
For example $w(19)=2$, since the symmetries of $F_{19}$
has two generators (the rank, or ``local'' dimension,
of the $\B{Z}$-module is two).
Alternatively, $p-1=19-1=18$ is a function supported on two points: $2$ and $3$.
The ``odd'' even prime $2$ is an initial element in this POSet, and plays the
role of a unit, as it will explained later on.
For now it will be treated on equal footing with the other primes,
even though the $2^k$ symmetry/factor of $p-1$ is ``special'', 
clearly a reflection linked to the orientation of the cycles of $Aut(F_p,+)$.
\begin{defin}
The {\em total weight} of a prime $p$ is 
$$W(p)=\sum_{p'|p-1} \nu_{p'}(p-1)).$$
\end{defin}
Here $\nu_p(n)=k(p)$ (p-adic valuation) is another notation for 
the value of the exponent function $k(p)$.

For example $W(19)=1+2=3$, since $19-1=2^1\cdot 3^2$.

These two functions will have natural interpretations when associating 
rooted trees to primes, reflecting the complexity of the primes, 
and of the structure of the corresponding primary fields.

We mention a few aspects to be studied, of possible interest to the reader:
a) refining the properties of arithmetic function, when taking into account the POSet structure of $P$;
b) the relation between Proth Theorem and
the structure of this POSet; c) the homology of the POSet of primes 
and its relation with the structure of rooted trees.

\section{Prime numbers and rooted trees}
To each prime $p$ we associate a rooted tree $t(p)$ representing
the hierarchy of generators of $F_p$, when applying $Aut_{Ab}$
repeatedly. 

This is essentially the Pratt tree, also called {\em prime certificate} \cite{Pratt},
without the explicit primitive root, or ``witness'', which is associated with its nodes.

We will define this recursively, 
at the level of prime factors of $Aut_{Ab}(F_p,+)$:
$$t(2)=\bullet, \quad t(p)=B^+(t(p_1),...,t(p_n)),\quad p-1=\prod_1^n p_i,$$
where $B^+$ is the operator which adjoins a common root to the input rooted trees.
For example
$$t(181)=B^+(\bullet, t(2), t(2), t(3), t(3), t(5)), \ t(3)=B^+(\bullet), 
\ t(5)=B^+(\bullet, \bullet).
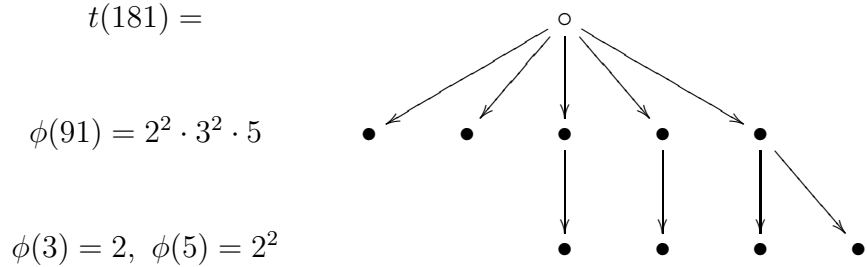
$$

$$\diagram
t(181)= & & & \circ \dllto \dlto \dto \drto \drrto & \\
\phi(91)=2^2\cdot 3^2\cdot 5 
  & \bullet & \bullet & \bullet \dto & \bullet \dto & \bullet \dto \drto \\
\phi(3)=2, \ \phi(5)=2^2 & & & \bullet & \bullet & \bullet & \bullet
\enddiagram$$

\begin{rem}
Note that multiplicity of a prime factor, e.g. $3^2$ in the above example,
is recorded as a multiplicity of the corresponding leaf.

Alternatively, to better reflect the underlying algebraic structure,
i.e. the elementary factor $\B{Z}_{3^2}$, one may chose to 
repeat the factor ``vertically'', 
attaching a {\em bamboo} of length equal to the multiplicity $k$, here $k=2$.
Then besides adjoining the leafs to a root, one needs to substitute
the trees corresponding to the elementary factors, to the nodes
of the bamboo (operation denoted $\circ$):
$$t(p)=B^+(b(k_1)\circ t(p_1), ..., b(k_l)\circ t(p_l)), 
\quad p-1=\prod_1^l p_j^{k_j}\ (p_i\ distinct).$$
In this case, a repetition ``scheme'' similar to Gerstenhaber composition, 
needs to be prescribed in order to define $\circ$.

At this early stages of development of the theory,
we chose the simpler definition above, using a ``horizontal'' scheme 
to record multiplicity.
\end{rem}

Once the mapping $t:P\to RT$, from prime numbers to rooted trees is defined,
it is natural to extend it linearly to a map between the 
integer lattice of prime numbers and 
the Hopf algebra of rooted trees $t:g_P\to H_{rt}$:
$$t(\sum_{p\in P} k_pX_p)=\sum_{k\in P} k_p t(p),$$
where the function $k:P\to \B{N}$ is finitely supported.

\subsection{An external composition on prime numbers}
The ``Euclid's trick'' (take products of primes, including 2, and add one),
almost allows to define a composition of primes.
It is the analog of root adjoining operator $B^+$ on rooted trees (and forests).

If $p$ and $q$ are primes, the product $(p-1)(q-1)$, 
which corresponds to joining their symmetries at categorical level,
might not correspond to the symmetries of a finite field, since
$$p*q=(p-1)(q-1)+1$$
is not necessarily a prime number, e.g. $3*5=9$.

Note, thought, that it defines an associative operation on $\B{Q}$,
and transport this definition on $g_P$, as a {\em fusion rule}
\footnote{$Hom_{Ab}$ is an internal functor.}:
$$X_p\star X_q=\sum_{r\in P} k(r)X_r, 
\quad (p-1)(q-1)+1=\prod_{r\in P}r^{k{r}}.$$ 
and extend it to a linear binary operation, still denoted $*$.
\begin{lem}
$*$ an associative binary operation on $g_P$ with identity
$X_2$.
\end{lem}
\begin{pf}
Associativity comes from the correspondence 
with the multiplicative group of rational numbers to be introduced in shortly,
and the fact that $*$ is an associative operation.

The identity is clearly $X_2$, since $2*q=(2-1)\cdot (q-1)+1=q$.
\end{pf}
\subsection{The rationals as a formal group}
To make the multiplicative group of rational numbers a formal group, 
``change coordinates'', defining the 
(multiplicative) ``position vector'' of $x\in Q$ by
$r_x=x-1$ (difference between $x$ and the neutral element).
Then define
$$x\star y=F(r_x,r_y)=r_x+r_y+r_x\cdot r_y.$$
\begin{prop}
There is a natural isomorphism of formal groups:
$$Exp:g_P\to (Q_+^\times, \cdot), \ Exp(X_p)=p-1,$$
between the additive formal group $(g_P,+)$ with addition:
$$F(X_p,X_q)=X_p+X_q,$$
and the multiplicative formal group $(\B{Q}_+^\times, \star)$
\cite{Wiki-formalgroups}.
The inverse is called the formal logarithm.
\end{prop}
\begin{rem}
The exponential implements a ``tangent map'' point of view, 
conform Lie Theory,
to reflect the action of the automorphism functor $Aut_{Ab}$,
and the fact that the basic finite fields $F_p$
are tangent spaces of the p-adic k-th order truncations of the 
p-adic ``quantum abelian groups''.
\end{rem}

\begin{rem}
This leads to a framework reminiscent of star products and deformation quantization:
$g_P$ plays the role of a Lie (bi)algebra, $Q$ is the associated group,
and p-adic numbers $Q_p=F_p((h))$ are a deformation of the laurent series 
in the formal parameter $h$.
This idea will be explored elsewhere.
\end{rem}

\subsection{A probabilistic framework}
The presence of this {\em two} fundamental ``God given'' structures,
prime numbers and rooted trees, and the link between them 
is highly intriguing; why two?

Taking as a working hypothesis that they are ``assymptotically'' 
the same in some sense, one could ponder whether prime numbers correspond 
to the cycles of the cohomology of rooted trees:
see the ``Universal cohomological problem \cite{ENCG-universalproblem},
or that rooted trees, as a basis in the Hopf algebra of trees,
could be a more fundamental probabilistic space to study primes
as a random variable, then the set of integers:
$$C:H_{rt}\to \B{Q}, \quad X:(PRT, ``\rho'')\to (\B{N}, card).$$
Here $PRT$ is the subspace in $g_P$ generated by {\em prime rooted trees},
the image of $t:P\to RT$, so that the image of the random variable $X$ 
is the set of primes $P\subset N$.
 
The linear mapping $C$ defined on the basis of rooted trees by
extending the random variable $X$ from $PRT$ to $RT$
should be a section of $t$.

The problem is to construct it using somehow the 
contraction 
$$\tau:g_P\to g_P, \quad \tau(X_p)=Div(Exp(X_p)),$$
where 
$$Div(n)=\sum k(p) X_p, \quad n=\prod_{p\in P} p^{k(p)},$$
is the {\em principal divisor} of the ``rational function''
$k:P\to \B{Z}$ associated to the natural number $n$
(more generally, to any non-zero rational number).

\section{On the Prime Number Theorem}
The general idea is to reduce the PNT, via linearization, 
to a ``simple'', Riemann-Roch-like statement,
easely prooved at the level of Hopf algebras of rooted trees,
its basis.

The Prime Number Theorem states that the probability density 
function of prime numbers is asymptotically $1/\ln n$:
$$P(1<X<n)=\pi(n)/n \sim 1/\ln(n), \quad n \in \B{N}.$$
Now $\ln(n)$ naturally ``lives'' in $g_P$,
the linearization of $\B{Q}$.
\begin{rem}
$\ln(n)=\sum k(p) \ln(p)$ is a ``divisor'' limiting
the number of basic elements $X_p$ in $g_p$,
in a way similar to Riemann-Roch Theorem.
This cannot be a simple ``coincidence''.
\end{rem}
The weight $\ln(p)$ present in the various formulas used in the proof
of PNT via the Riemann zeta function is in fact associated with
a measure.
It appears under the disguise of the {\em Mangold ``function''}
$$\Lambda(p^k)=\ln(p), \quad (or \Lambda(n)=0\ otherwise).$$
Recall the other Chebyshev function:
$$\Psi(x)=\sum_{n<x} \Lambda(n),$$ 
and note that it can be interpreted as the integral on $(g_P, \ln(p) dp)$,
of a cutoff exponent function $K(p)=\lceil x/p \rceil$.

As for the ``definite integrals'' (the above dependence on $x$ can be 
interpreted as an anti-derivative), the following diagram relates
the functional model of the rational numbers with the
usual context:
$$\diagram
(g_P, \ln(p)dp) \drto^{\int} \rto^{exp} & (\B{Q}_+,\cdot) \dto^{\ln(r)} \\
& (\B{R}, +, dx)
\enddiagram$$
$$p=exp(X_p), \quad \int_P k=\sum k(p)\ln(p)=\ln(n),$$
where the divisors $k(p)$ are interpreted as functions,
and the exponential is the natural isomorphism determined by the
bijection between the two basis of the two free abelian groups.

The natural measure on $g_P$ ``close'' to the 
measure $\ln(p)dp$ is the ``size'' of the corresponding 
rooted tree $t_p$, i.e. its degree:
$$deg(t_p)=number \ of \ nodes.$$
\begin{conj}
$$\ln(p) \sim deg(t_p).$$
\end{conj}
\begin{pf}
A preliminary experimental check involving students was performed as part of the
ISU Summer Research Academy 2012-2013 \cite{ISU-SRA}. 
\end{pf}
The proportionality factor obtained could be due to a Gaussian distribution 
of the exponents of primes as symmetries, which were not accounted for in this
particular numerical exploration.

\section{Conclusions}
Some new ideas are introduced in order to study the prime numbers.

First, categorification teaches us to study the objects themselves,
here the basic finite fields, not their iso classes; and here natural numbers 
cannot even distinguish between $Z_{p^k}$ and $Z_p^k$.

From the start, studying the symmetries of simple abelian groups leads to a
partial order on the set of prime numbers.
This in turn, allows for a correspondence with rooted trees and 
several grading functions, with a corresponding concept of complexity of 
prime numbers.

Next, viewing the multiplicative group of rational numbers as functions on the set of primes, leads to 
alternative approaches for a proof of the Prime Number Theorem,
where the POSet structure allows to structure the internal ``production'' of new
primes out of old.
Whether this will lead to a direct \underline{and} conceptual proof of the PNT,
remains to be seen.

Finally, the Lie-Klein geometric interpretation, which is suited for generalization
in the context of formal groups and deformation theory,
leads to consider algebraic quantum groups as the proper framework for studying
adeles and the duality between primes and zeta zeroes,
which is ``behind'' the Riemann Hypothesis.

Overall, a modern approach in Multiplicative Number Theory is needed, beyond
the traditional analytic approach, which has ``too much analysis'';
and since quantum theory is envisioned to help prove the Riemann Hypothesis,
deformation theory and path integrals should be used as main tools.

\section{Acknowledgments}
I thank Prof. Carl Pomerance for his comments, including the relation with Pratt trees.


\end{document}